\documentclass[11pt]{article}

\usepackage{amsmath,amssymb}
\usepackage{graphicx}
\usepackage{xcolor}
\usepackage{xspace}
\usepackage{natbib}
\usepackage[nameinlink]{cleveref}
\usepackage{pgffor}

\usepackage[
left=2.5cm,
right=2.5cm,
top=2.5cm,
bottom=2.5cm
]{geometry}

\newtheorem{dfn}{Definition}

\newtheorem{prp}{Proposition}

\newtheorem{exm}{Example}
\newtheorem{rem}{Remark}
\newtheorem{thm}{Theorem}

\newenvironment{prf}[1][\unskip]{%
	\textbf{Proof #1.~}}{\hfill$\blacksquare$}

\newcommand{\almostsurely}{\text{\xspace a.\,s.}\xspace}				
\newcommand{\median}[2][{}]{\text{Med}_{#1}\left[#2\right]}				

\newcommand{\Noise}{B}													
\newcommand{\noise}{b}													

\DeclareMathOperator*{\esssup}{ess\,sup}								
\newcommand{\goaldist}[1][{}]{d_{\G#1}}									

\newcommand{\mquad}[1]{%
	\foreach \n in {1,...,#1} {%
		\qquad%
	}%
}

\makeatletter
\newcommand{\pushright}[1]{\ifmeasuring@#1\else\omit\hfill$\displaystyle#1$\fi\ignorespaces}
\newcommand{\pushleft}[1]{\ifmeasuring@#1\else\omit$\displaystyle#1$\hfill\fi\ignorespaces}
\makeatother

\newcommand{\nomitem}[2]{#1 & #2 \\ \hline}

\newcommand{\diff}{\mathop{}\!\mathrm{d}}								
\newcommand{\eps}{{\varepsilon}}										
\newcommand{\ball}{{\mathcal B}}										
\newcommand{\diam}{{\text{diam}}}										
\newcommand{\image}{{\text{Image}}}										
\newcommand{\inv}{\ensuremath{^{-1}}}									
\newcommand{\nrm}[1]{\left\lVert#1\right\rVert}							
\newcommand{\abs}[1]{\left\lvert#1\right\rvert}							
\newcommand{\floor}[1]{\left\lfloor#1\right\rfloor}						
\newcommand{\cl}[1]{\text{cl}\left({#1}\right)}							
\newcommand{\E}[2][{}]{\mathbb E_{#1}\left[#2\right]}					
\newcommand{\PP}[1]{\mathbb P\left[#1\right]}							
\newcommand{\low}{{\text{low}}}											
\newcommand{\up}{{\text{up}}}											
\newcommand{\ra}{\rightarrow}											
\newcommand{\la}{\leftarrow}											
\newcommand{\ie}{\unskip, i.\,e.,\xspace}								
\newcommand{\eg}{\unskip, e.\,g.,\xspace}								
\newcommand{\sut}{\text{s.\,t.\,}}										
\newcommand{\wrt}{w.\,r.\,t. \xspace}									

\let\emptyset\varnothing
\newcommand{\N}{{\mathbb{N}}}											
\newcommand{\Z}{{\mathbb{Z}}}											
\newcommand{\R}{{\mathbb{R}}}											
\newcommand{\state}{s}													
\newcommand{\State}{S}													
\newcommand{\states}{\mathbb S}											
\newcommand{\action}{a}													
\newcommand{\Action}{A}													
\newcommand{\actions}{\mathbb A}										
\newcommand{\traj}{z}													
\newcommand{\Traj}{Z}													
\newcommand{\policy}{\pi}												
\newcommand{\policies}{\Pi}												
\newcommand{\transit}{p}												
\newcommand{\cost}{c}													
\newcommand{\G}{\ensuremath{\mathbb{G}}}								
\newcommand{\K}{\ensuremath{\mathcal{K}}\xspace}						
\newcommand{\KL}{\ensuremath{\mathcal{KL}}\xspace}						
\newcommand{\Kinf}{\ensuremath{\mathcal{K}_{\infty}}\xspace}			
\newcommand{\T}{\mathcal T}												
\newcommand{\spc}{{\,\,}}												

\begin{document}

\title{Some Remarks on Stochastic Converse Lyapunov Theorems}

\author{
	Pavel Osinenko\thanks{Authors contributed equally. Corresponding author: Pavel Osinenko, e-mail: \texttt{p.osinenko@gmail.com}}
	\and
	Grigory Yaremenko\footnotemark[1]
}

\date{} 

\maketitle

\begin{abstract}
	In this brief note, we investigate some constructions of Lyapunov functions for stochastic discrete-time stabilizable dynamical systems, in other words, controlled Markov chains.
The main question here is whether a Lyapunov function in some statistical sense exists if the respective controlled Markov chain admits a stabilizing policy.
We demonstrate some constructions extending on the classical results for deterministic systems.
Some limitations of the constructed Lyapunov functions for stabilization are discussed, particularly for stabilization in mean.
Although results for deterministic systems are well known, the stochastic case was addressed in less detail, which the current paper remarks on.
A distinguishable feature of this work is the study of stabilizers that possess computationally tractable convergence certificates.
\end{abstract}




%
\label{sec_notation}
\section*{Notation}

\renewcommand{\arraystretch}{1.2}
\setlength{\tabcolsep}{10pt}

\begin{center}
	\begin{tabular}{|l|p{0.72\linewidth}|}
		\hline
		\textbf{Symbol} & \textbf{Meaning} \\ \hline
		
		\nomitem{$\cl{\bullet}$}{Closure of $\bullet$}
		\nomitem{$\diam_\bullet$}{Diameter of $\bullet$, i.e. $\sup_{x, y \in \bullet} \nrm{x - y}$}
		\nomitem{$\ball_v(x)$}{Closed ball of radius $r$ centered at $x$; ``$(x)$'' is omitted when $x = 0$}
		\nomitem{$\K$}{Set of all continuous, strictly increasing, positive-definite functions over positive reals}
		\nomitem{$\Kinf$}{Set of all unbouned $\K$-functions}
		\nomitem{$\bullet \sim \bullet$}{$\bullet$ (left) is distributed according to $\bullet$ (right)}
		\nomitem{$\bullet := \bullet$}{$\bullet$ (left) is set as $\bullet$ (right)}
		\nomitem{$\floor{\bullet}$}{Floor of $\bullet$}
		\nomitem{\almostsurely}{Almost surely}
		\nomitem{$\median{\bullet}$}{Median of $\bullet$}
		\nomitem{$\neg \bullet$}{Complement of $\bullet$}
		\nomitem{$d_{\bullet}(\bullet)$}{Distance from set $\bullet$ (subscript) to point $\bullet$ (parentheses)}
		\nomitem{$\states$}{State space}
		\nomitem{$\actions$}{Action space}
		\nomitem{$p$}{Transition probability}
		\nomitem{$\KL$}{See \Cref{dfn_kappaell}}
		\nomitem{$\mathbb{G}$}{Goal set}
		\nomitem{$\Omega$}{Sample space}
		\nomitem{$\Sigma$}{Event space}
		\nomitem{$\mathbb{P}[\bullet]$}{Probability of $\bullet$}
		\nomitem{$\mathbb{E}[\bullet]$}{Expected value of $\bullet$}
		\nomitem{$\mathcal{E}_0^\bullet$}{Vanishing locus of $\bullet$}
		\nomitem{$\image(\bullet)$}{Image of map $\bullet$}
		\nomitem{$\bullet_{0:N}$}{$\bullet_0, \bullet_1, ..., \bullet_N$ (analogously for superscript)}
		
	\end{tabular}
\end{center}


\section{Introduction}
\label{sec_relwork}

Converse Lyapunov theorems for deterministic dynamical systems have been extensively studied in the literature, with notable contributions including the seminal works of \cite{Lin1996smoothconverse,Jiang2002converseLyapun,Kellett2004WeakconverseLa}. While corresponding results for stochastic systems are also well-established and continue to be developed, they remain relatively less familiar compared to their deterministic counterparts. Nevertheless, the significance of converse Lyapunov theorems extends beyond control theory, finding important applications in data science, particularly in stochastic approximation (\cite{Vidyasagar2022newconverseLy}).

For continuous stochastic dynamical systems, foundational results can be found in the classical works by \cite{Cesaroni2006converseLyapun,Liu2011converseLyapun}, while more recent contributions include those of \cite{Rajpurohit2016Lyapunovconver,Pepe2014Directconverse,Majumdar2024NecessarySuffi}. As is typical in the analysis of stochastic differential equations, these works employ growth conditions as key technical assumptions. Similarly, \cite{Pepe2014Directconverse} utilized growth conditions to establish converse Lyapunov function results for functional difference systems.

In the context of controlled Markov chains, particularly significant contributions have been made by \cite{Teel2014converseLyapun,Teel2012converseLyapun}. These works address almost sure global asymptotic recurrence through sophisticated recursive constructions. The present work adopts the more straightforward methodology developed in \cite{Jiang2002converseLyapun}, utilizing the $\kappa$-function framework of \cite{Khalil1996NonlinearSyste,Sontag1998Commentsintegr}. 

\cite{Subbaraman2013converseLyapun} investigated Lyapunov function constructions for strong global recurrence, where solutions return almost surely to an open and bounded set infinitely often. Their approach employed convexification of bounding $\kappa$-functions to leverage Jensen's inequality effectively. Meanwhile, \cite{Haddad2022Lyapunovtheore,Haddad2020StochasticSemi} examined systems with additive diffusion that vanishes at the equilibrium.

Of particular interest are constructive approaches utilizing linear or mixed-integer programming (\cite{Baier2012Linearprogramm,Hafstein2005constructiveco}), although these methods have thus far been limited to the deterministic case.

Finally, it is worth noting that converse results have also been established for stochastic differential inclusions \cite{Teel2012Stochasticdiff}.

As compared to existing discrete stochastic converse Lyapunov results the presented converse theorems are distinguished by
the following:
\begin{itemize}
\item The results concern $\eta$-improbable stabilization, i.e. policies that may fail to stabilize the system with a certain probability. This makes the results of the paper applicable to a significantly more general class of converse problems, like those that involve not only stabilizers, but also policies, whose stabilizing properties can be broken by the stochastic phenomena in the system's dynamics.
\item The assumptions on system dynamics in the current work are weaker than those found in the literature \eg \cite{Subbaraman2013converseLyapun} (for instance, local boundedness of transitions is not required).
\item The results are applicable not only to stabilization of a single equlibrium, but also to that of arbitrary attractors.
\end{itemize}

\section{Background and problem statement}
\label{sec_problem}


Consider the following controlled Markov chain:
\begin{equation}
	\label{eqn_mdp}
	\left(\states, \actions, \transit \right),
\end{equation}
where:
$\states$ is the \textit{state space}, assumed as a normed vector space of all states of the given environment;
$\actions$ is the \textit{action space}, that is a set of all actions available to the agent, it may be discrete or continuous, compact or non-compact;
$\transit : \states \times \actions \times \states \ \rightarrow \ \R$ is the \textit{transition probability density function} of the environment, that is such a function that $\transit(\bullet \mid \state_{t}, \action_{t})$ is the probability density of the state $s_{t + 1}$ at step $t+1$ conditioned on the current state $\state_{t}$ and current action $\action_{t}$;
%
Let $(\Omega, \Sigma, \PP{\bullet})$ be the probability space underlying \eqref{eqn_mdp}, and $\E{\bullet}$ be the respective expected value operator.
The problem of stabilization is to find a policy $\policy_0$ from some space of admissible policies $\policies_0$ that drives the state into some goal set $\G \subset \states$ in some suitable statistical sense.
Let the goal contain the origin for simplicity \ie $0 \in \G$ and the distance to it be denoted $\goaldist(\state) := \inf\limits_{\state' \in \G} \nrm{\state - \state'}$.
%
%
The stabilization problem may be concerned with a designated initial state $\state$, a distribution thereof, or even the whole state space.
A policy may be taken as a probability density function or as an ordinary function (cf. Markov policy).
The said statistical sense may be understood differently depending on context.
For instance, one may seek stabilization with probability not less than $1-\eta, \eta>0$, in median etc.
We call a policy $\policy_0 \in \policies_0$ a \textit{$\G$-stabilizer} if setting  $\Action_t \gets \policy_0(\bullet \mid \state_t), \forall t$ implies that the distance between $\State_t$ and $\G$ tends to zero over time in a statistical sense.
Formally, we use the following definition.

\begin{dfn}
	\label{dfn_stabilizer}
	A policy $\policy_0 \in \policies_0$ is called a \textit{uniform $\eta$-improbable $\G$-stabilizer}, $\eta \in [0, 1)$, or simply a \textit{stabilizer} (if the goal set is implied by the context) if
	\begin{equation}
		\label{eqn_goallim}
		 \PP{\lim_{t \ra \infty} \goaldist(\State_t) = 0 \mid \Action_t \sim \policy_0(\bullet \mid \State_t)} \ge 1 -\eta \ \forall \state_0 \in \states,
	\end{equation}
	where the limit in \eqref{eqn_goallim} is compact-convergent \wrt $\state_0$ and
	\begin{equation}
		\label{eqn_goalstab}
		\begin{aligned}
			& \exists \eps_0 \ge 0 \spc \forall t \ge 0 \spc \Action_t \sim \policy_0(\bullet \mid \State_t) \implies \\
			& \mquad{2} \forall \eps \ge \eps_0 \spc \exists \delta \ge 0 \spc \goaldist(\state_0) \le \delta \\
			& \mquad{3} \implies  \PP{\goaldist(\State_t) \le \eps} \ge 1 -\eta \ \forall t \geq 0, 
		\end{aligned}
	\end{equation}	
	where $\delta$ is unbounded as a function of $\eps$ and $\delta = 0 \iff \eps = \eps_0$.
	The presence of an $\eps_0$ means we do not insist on $\G$ being invariant.
	Thus, this extra condition only means that the state overshoot may be uniformly bounded over compacts in $\states$.
\end{dfn}



The above can be viewed as a generalization of global weak asymptotic stability in probability (\cite{Khasminskii2011StochasticStab}). This definition is distinguished by its applicability to the case of set-attractors and stochastic phenomena that have a non-zero probability of preventing stabilization.

Statistical confidence bounds by \eg Hoeffding’s inequality in terms of the policy runs on the system\footnote{Cf. \cite{Hertneck2018Learningapprox}} can be used to verify \eqref{eqn_goallim} or \eqref{eqn_meanstabilizer}.
The \textbf{aim of this work} is to investigate simple constructions of Lyapunov functions in the spirit of \cite{Jiang2002converseLyapun} for stabilized controlled Markov chains.
Given that a policy $\policy_0$ exists that solves a stabilization problem of \eqref{eqn_mdp}, we would like to find a Lyapunov function in the respective statistical sense.
Furthermore, we would like to remark on stabilization by means of a constructed Lyapunov function.
The special focus here is on $\eta$-improbable stabilizers which possess special tractability properties.
This is what primarily distinguishes the current work.
The details follow in \Cref{sec_thms}.

\section{Results}
\label{sec_thms}


\begin{dfn}
	\label{dfn_kappaell}
	A continuous function $\beta : \R_{\ge0} \times \R_{\ge0} \ra \R_{\ge0}$ is said to belong to space $\KL$ if
	\begin{itemize}
		\item $\beta(v, \tau)$ is a $\K$ function with respect to $v$ for each fixed $\tau$. 
		\item $\beta(v, \tau)$ is a strictly decreasing function with respect to $\tau$ for each fixed $v$.
		\item For each fixed $v$ it holds that $\lim\limits_{\tau \ra \infty}\beta(v, \tau) = 0$.
	\end{itemize} 
\end{dfn}

The following technical result will be needed later on. 
It may be considered as a generalization of the classical results on existence of $\KL$ convergence certificates for uniformly globally stable systems (\cite{Lin1996smoothconverse}). 
Unlike the latter, \Cref{prp_lim2kappa} is system-agnostic.

\begin{prp}
	\label{prp_lim2kappa}
	Consider a function $\Phi: \R^n \times \R_{\ge 0} \ra \R_{\ge 0}$ with the property that
	\begin{equation}
		\label{eqn_genericlim}
		\begin{aligned}
            & \forall \text{ compact } \states_0, \eps>0 \spc \exists \T>0 \ \forall \tau \ge \T\!, \state \in \states_0
			\ \Phi(\state, \tau) \le \eps,  \\  
			& \forall \eps \ge 0 \spc \exists \delta \geq 0 \spc \nrm{\state} \le \delta \implies \forall \tau \ge 0 \spc \Phi(\state, \tau) \le \eps.
		\end{aligned}		
	\end{equation}
	where $\delta$ is unbounded as a function of $\eps$. 
%
	Then, there exists a $\KL$ function $\beta$ with the following property:
	\begin{equation}
		\label{eqn_exitskl}
		\forall \state \in \R^n, \tau \ge 0 \spc \Phi(\state, \tau) \le \beta(\nrm{\state}+C_0, \tau),
	\end{equation}
	where $C_0 \ge 0$ is a constant that equals zero if $\delta$ is positive outside zero as a function of $\eps$.
\end{prp}

\begin{prf}
	First, we may assume that $\delta$ is a non-decreasing function of $\eps$.
	To see this, fix $0 < \eps_1 \le \eps_2$.
	Then, 
	\[
		\nrm{\state} \le \delta(\eps_1) \implies \forall \tau \ge 0 \spc \Phi(\state, \tau) \le \eps_1,
	\]
	from which it trivially follows that $\Phi(\state, \tau) \le \eps_2$ also.
	Hence, redefining $\delta(\eps) := \sup\limits_{\eps' \in [0, \eps]}\delta_\text{old}(\eps')$ is valid without loss of generality.	
	Second, $\T_{\states_0}$, for a fixed compact set $\states_0$ is a non-increasing function of $\eps$.
	Again, fix some $0 < \eps_1 \le \eps_2$ and observe that, for $\state$ in $\states_0$, $\tau \ge \T_{\states_0}(\eps_1) \implies \Phi(\state, \tau) \le \eps_1$
	implies trivially $\Phi(\state, \tau) \le \eps_2$ whenever $\tau \ge \T_{\states_0}(\eps_2)$.
	Now, let us consider $\T$ as a function of two arguments -- radius $v$ of the minimal ball containing the given compact set $\states_0$ and $\eps$ -- for convenience.	
	Then, it is a non-increasing function in the first argument.
	Yet again, fix any $0 < \eps_1 \le \eps_2$ and $v >0$.
	Take two compact sets $\states_0, \states_0'$ with $\states_0 \subseteq \states_0'$ and $\states_0' \subseteq \ball_v, v>0$.
	Then, for any $\eps>0$, we have
		$\state \in \states_0', \tau \ge \T(v, \eps) \implies \Phi(\state, \tau) \le \eps$
	and, simultaneously,
		$\state \in \states_0, \tau \ge \T(v, \eps) \implies \Phi(\state, \tau) \le \eps.$
	Furthermore, we redefine $\delta$ to be continuous, strictly increasing and positive on $[v_0, \infty)$, where $v_0 := \sup \mathcal E^\delta_0$ and $\mathcal E^\delta_0:= \{ \eps : \delta(\eps) = 0 \}$.
	To see this, observe that $\delta$ is positive everywhere outside zero, non-decreasing and bounded by $\delta(b)$ on any $[a, b]$.
	Hence, $\delta$ is Riemann integrable on an arbitrary bounded interval precisely because the sets of discontinuities of non-decreasing functions have measure zero on compact domains.
	Consider a function $\hat \delta$ defined by
	\[
		\hat \delta(\eps) := \frac{1}{\eps + 1}\int_{0}^{\eps} \delta(\eps') \diff \eps'.
	\]
	Evidently, $\hat\delta$ is continuous. 
	Notice, that $\hat\delta(\eps)$ is zero whenever $\delta(\eps)$ is zero and for $\eps - v_{0} > h > 0$ it holds due to $\delta$ being non-decreasing that
	\begin{equation}
		\begin{aligned}
			\hat\delta(\eps) = \frac{1}{\eps + 1}\int_{0}^{\eps} \delta(\eps') \diff \eps' \geq \frac{1}{\eps + 1}\int_{\eps - h}^{\eps} \delta(\eps') \geq \\
			\frac{1}{\eps + 1}h\delta(\eps - h) > 0.
		\end{aligned}
	\end{equation}
	Therefore $\hat\delta(\eps)$ is positive whenever $\delta(\eps)$ is positive except possibly at $v_0$. 
	Now, note that assuming $\eps > 0$ by mean value theorem we have
	\begin{equation}
		\begin{aligned}	
			 & \hat\delta(\eps) = \frac{1}{\eps + 1}\int_{0}^{\eps} \delta(\eps') \diff \eps' \leq \frac{1}{\eps}\int_{0}^{\eps} \delta(\eps') \diff \eps'= \delta(\hat\eps) \\
			 & \pushright{\leq \delta(\eps), \text{ where }\hat\eps \in [0, \eps]}.
		\end{aligned}	 
	\end{equation}
	Since $\delta(0) = \hat\delta(0)$, it is true that $\hat\delta$ bounds $\delta$ from below.
	Also, for $\eps \geq 1$ by mean value theorem

	\begin{equation}
		\begin{aligned}	
			& \hat\delta(\eps) = \frac{1}{\eps + 1}\int_{0}^{\eps} \delta(\eps')\diff \eps' \geq \frac{1}{2\eps}\int_{0}^{\eps} \delta(\eps') \diff \eps'\geq \\
			& \mquad{5} \frac{1}{4} \frac{1}{\frac{1}{2}\eps}\int_{\frac{\eps}{2}}^{\eps} \delta(\eps')\diff \eps' \geq \frac{1}{4}\delta\left(\frac{\eps}{2}\right),
		\end{aligned}	
	\end{equation}

	therefore $\hat\delta$ is unbounded.
	Finally, for $h > 0$, $\eps \notin \mathcal E^\delta_0$
	\begin{align*}
		& h\int_{0}^{\eps}\delta(\eps')\diff \eps'\leq h\eps\delta(\eps) < h(\eps + 1)\delta(\eps) \\
		& \mquad{4} \leq (\eps + 1)\int_{\eps}^{\eps + h}\delta(\eps') \diff \eps' \\
		& \implies (\eps + 1)\int_{0}^{\eps}\delta(\eps')d\eps' + h\int_{0}^{\eps}\delta(\eps')\diff \eps' < \\
		& \mquad{2} (\eps + 1)\int_{0}^{\eps}\delta(\eps')\diff \eps' + (\eps + 1)\int_{\eps}^{\eps + h}\delta(\eps') \diff \eps'\\
		& \implies
	\end{align*}
	\begin{align*}		
		& (\eps + 1 + h)\int_{0}^{\eps}\delta(\eps')\diff \eps' < (\eps + 1)\int_{0}^{\eps + 1 + h}\delta(\eps')\diff \eps'\\
		& \implies \\
		& (\eps + 1)(\eps + 1 + h)\hat\delta(\eps) < (\eps + 1)(\eps + 1 + h)\hat\delta(\eps + h)\\
		& \implies \hat\delta(\eps) < \hat\delta(\eps + h),
	\end{align*}
	therefore $\hat\delta$ is strictly increasing outside of $\cl{\mathcal E^\delta_0}$.	
	
	Let us now fix an $\state \in \R^n, \state \ne 0$. 
	Then, by the compact convergence, it holds that
	\[
		\forall \state' \in \ball_{\nrm{\state}}, \eps > 0, \tau \ge \T(\nrm{\state}, \eps) \spc \Phi(\state', \tau) \le \eps.
	\]
	Let $\eps_x := \delta\inv(\nrm{\state})$ which is possible since $\delta$ is now assumed continuous and strictly increasing outside $\cl{\mathcal E^\delta_0}$.
	Notice $\eps_x = v_0$ if $\state=0$ by the above redefinition of $\delta$.
	In particular, it is zero if the original $\delta$ was positive outside zero.
	From the overshoot bound condition, second line of \eqref{eqn_genericlim}, it holds that:
	\[
		\forall \state' \in \ball_{\nrm{\state}}, \tau \ge 0 \spc \Phi(\state', \tau) \le \eps_x.
	\]
	Therefore, it is safe to assume that $\T(\nrm{\state}, \eps)$, as a function of $\eps$, has a compact support since $\eps_x$ is never exceeded by $\Phi$ as shown above.
	
	Now, define inductively the numbers $\eps_j, j \in \N$ starting with $\eps_1 := \eps_x$:

	\[
		\eps_j :=\! 
		\begin{cases}
			\sup \mathcal E^\T_{\state,j}, \text{ if } \mathcal E^\T_{\state,j}\! :=\! \{ \eps\! : \T(\nrm{\state}, \eps)\! \in\! [j-1,j] \}\! \ne\! \emptyset, \\
			\eps_{j-1}, \text{ otherwise}.
		\end{cases}
	\]

	By the convergence property \eqref{eqn_genericlim}, it holds that \\ 
	$\limsup\limits_{\eps \ra \infty} \T(\nrm{\state}, \eps) = 0$ and $\T$ is defined for any $\eps$.
	Hence, $\T$ may not ``blow up'' into infinity at any $\eps>0$ and so the above construction of the numbers $\eps_j$ is valid.

	Now, let $\beta'_v, v = \nrm{\state}$ be a function, such that its graph is a polygonal chain specified by \\ $(0, \eps_1), (1, \eps_1), (2, \eps_2), (3, \eps_3) \dots$ \ie
	\begin{align*}
		& \beta'_v(\tau) = \eps_{1} + \sum_{i = 1}^{\infty}(\eps_{i + 1} - \eps_{i})\chi(\tau - i), \\ 
		& \mquad{4} \text{ where } \chi(\tau) = \frac{1}{2}(1 + \lvert \tau \rvert - \lvert \tau - 1 \rvert).
	\end{align*}
	Evidently, $\beta'_v$ is continuous.
	A setting $\beta'_v \mapsto \beta'_v + C_1 v e^{- C_2 \tau}$, with $C_1, C_2 >0$ arbitrary, makes it also strictly decreasing.
	Notice that any other function $\T'(\nrm{\state}, \eps)$, defined for all $\eps > 0$, with $\T' \ge \T$ and $\limsup\limits_{\eps \ra \infty} \T'(\nrm{\state}, \eps) = 0$ is a valid certificate for the convergence property as per \eqref{eqn_genericlim}.	
	We thus see that the inverse of $\beta'_v$, which exists, is such a valid certificate. Now, let

	\begin{align}
		& \beta(v, \tau)\! :=\! (v \! - \! \floor{v})\beta_{\floor{v}\! +\! 2}(v, \tau)\! +\! (\floor{v}\! +\! 1\! -\! v)\beta_{\floor{v} + 1}(v, \tau), \notag \\ 
		& \text{where } \beta_i(v, \tau) := 2^i\max\limits_{j = 1}^i \beta^*_j(v, \tau), \notag \\
		& \beta^*_i(v, \tau) := \min(\beta'(i, \tau), \xi(v)), \notag\\
		& \xi(v) := \left\{\begin{array}{ll}
		                  v, \text{ if }v < \delta^{-1}(0),\\
		                  \delta^{-1}(v - \delta^{-1}(0)), \text{ otherwise.}
		                \end{array}\right.
	\end{align}

	It is true that $\beta'(i, \bullet)$ provides a valid upper bound on $\Phi(\state, \tau)$ not only for $\nrm{\state} = i$, but also for $\nrm{\state} < i$. 
	However the largest overshoot is no greater than $\xi(\nrm{\state} + \delta^{-1}(0))$, thus $\min(\beta'(i, \tau), \xi(\nrm{\state} + \delta^{-1}(0)))$ sharpens the former bound. 
	Evidently $\beta^*_{i}$ are bounded above by $\beta_{i}$ and therefore $\nrm{\state} \leq i$ implies $\beta_{i}(\nrm{\state} + \delta^{-1}(0), \tau) \geq \Phi(\state, \tau)$. 
	Finally, $\beta(\nrm{\state} + \delta^{-1}(0), \bullet)$ is a convex combination of $\beta_{\lfloor \nrm{\state} \rfloor + 1}(\nrm{\state} + \delta^{-1}(0), \bullet)$ and $\beta_{\lfloor \nrm{\state} \rfloor + 2}(\nrm{\state} + \delta^{-1}(0), \bullet)$, which are both upper bounds for $\Phi(\state, \tau)$. 
	The latter thus also applies to $\beta(\nrm{\state} + \delta^{-1}(0), \bullet)$.
	Note that $\beta_i$ are composed of continuous functions and therefore are continuous themselves. 
	At the same time, when $\beta$ is restricted to an arbitrary segment $[l - 1, l + 1], \ l \in \mathbb{N},$ it can be represented as $\chi_3(v - l + 1, \beta_{l}(v,\tau), \beta_{l + 1}(v,\tau), \beta_{l + 2}(v, \tau))$ where $\chi_3(t, k_{1}, k_{2}, k_{3}) = k_{1} + (k_{2} - k_{1})\chi(t - 1) + (k_{3} - k_{2})\chi(t - 2)$.
	Since $\chi$ and $\beta_{i}$ are continuous, then $\beta$ is continuous on an arbitrary segment $[l - 1, l + 1], \ l \in \mathbb{N},$ and thus $\beta$ is continuous on the entirety of its domain.
	Notice, that $\beta^*_{i}$ are non-decreasing w.r.t $v$ and thus so are $\beta_{i}$.
	Now, let $0 < v < v' < \lfloor v\rfloor + 1$ and observe that

	\begin{align}
		&\beta(v', \tau) - \beta(v, \tau) \geq \notag\\
		&(v'- v)\beta_{\lfloor v \rfloor + 2}(v, \tau) + (v - v' )\beta_{\lfloor v\rfloor + 1}(v, \tau) = \notag\\ 
		&(v' - v)(\beta_{\lfloor v \rfloor + 2}(v, \tau) - \beta_{\lfloor v\rfloor + 1}(v, \tau)) \geq \notag\\
		&(v' - v)(2\beta_{\lfloor v \rfloor + 1}(v, \tau) - \beta_{\lfloor v\rfloor + 1}(v, \tau)) = \notag\\ 
		&(v' - v)\beta_{\lfloor v \rfloor + 1}(v, \tau) \ge \notag\\
		&(v' - v)\min(\beta'(\lfloor v \rfloor + 1, \tau), \delta^{-1}(v)) > 0.
	\end{align}

	By continuity, this also extends to $v = 0$, therefore $\beta$ is increasing \wrt $v$.
	Also note that for $l \in \mathbb{N}$
	\begin{equation}
		\beta(l, \tau) = \beta_{l + 1}(l, \tau) \geq 2^{l + 1} \min(\beta'(2, \tau), \delta^{-1}(1)),
	\end{equation}
	thus $\beta$ is unbounded \wrt $v$ for every $\tau$.
	
	It is true that $v = 0$ if and only if $\beta(v, \tau) = 0$.
	Indeed,
	\begin{align*}
		& \beta(0, \tau) = 2 \min(\beta'(1, \tau), \delta^{-1}(0)) = 0; \\ 
		& v > 0 \implies \beta(v, \tau) \geq \beta'(\lfloor v\rfloor + 1, \tau) > 0  \ \lor \\
		& \pushright{\beta(v, \tau) \geq \delta^{-1}(v) > 0}.
	\end{align*}
	Trivially for $\alpha > 0, v > 0$, function $\alpha\beta_{i}(v, \bullet)$ is decreasing, positive and tending to $0$. 
	For whole values of $v$, it is true that $\beta(v, \bullet) = \beta_{v + 1}(v, \bullet)$, otherwise $\beta(v, \bullet)$ can be represented as $\alpha_{1}\beta_i(v, \bullet) + \alpha_{2}\beta_j(v, \bullet)$, thus $\beta(v, \tau)$ is decreasing, positive and tending to $0$ \wrt $\tau$ for every fixed $v > 0$ since the same applies to the terms.
	Thus it has been established that $\beta$ is a $\Kinf$ function \wrt $v$ for each $\tau$, while also being decreasing, positive and tending to $0$ w.r.t. $\tau$, which proves that $\beta$ is indeed a $\KL$-function. 
	At the same time it has been demonstrated that $\Phi(\state, \tau) \leq \beta(\nrm{\state} + C_0, \tau)$, where $C_0 = \delta^{-1}(0)$.
\end{prf}

\begin{rem}
	\label{rem_brokenlinedelta}
	The continuous and strictly increasing redefinition of $\delta$ by means of Riemann integrals admits an alternative construction as follows.
	Define for each interval $[k, k+1], k \in \Z_{\ge 0}$, the numbers:
	\[
		d_{k+1} := \inf_{\eps \in [k, k+1]} \cl{\{ \delta(\eps) \}}.
	\]
	This setting is possible due to the fact that $\delta$ is bounded for all $\eps$ and $\lim\limits_{\eps \ra \infty} \delta(\eps) = \infty$.
	That $\delta$ is bounded for any $\eps$ and tends to infinity as $\eps$ does follows from the trivial observation that $\delta \le \eps$ always and $\delta$ was assumed unbounded.
	Notice that any other unbounded function $\delta'$ with $\delta' \le \delta$ is a valid certificate for the overshoot bound property as per the second line of \eqref{eqn_genericlim}.
	Hence, we safely redefine $\delta$ to be such that its graph is a polygonal chain specified by $(0, d_1), (1, d_1), (2, d_2), (3, d_3) \dots$\ie
	\begin{align*}
		& \delta(\eps) := d_{1} + \sum_{i = 1}^{\infty}(d_{i + 1} - d_{i})\chi(t - i), \\
		& \mquad{3} \text{ where } \chi(\eps) := \frac{1}{2}(1 +  \abs{\eps}  -  \abs{\eps-1} ).
	\end{align*}
	Notice neither the original nor the thus redefined $\delta$ needs to be zero at zero and/or strictly increasing.
	However, we can make it strictly increasing on $[v_0, \infty)$ by defining $\delta$ to be
	\[
		\hat \delta(\eps) := \int_{0}^{\eps} \delta(\eps') \diff \eps'.
	\]
	similar to the construction in the proof, except for the factor, and defined only on the unit interval.
	To this end, let us assume, without loss of generality, that $\cl{\mathcal E^\delta_0}$ is well contained in the unit interval.
	Hence, unlike the polygonal chain construction above, we do not redefine $\delta$ to be the horizontal line segment on the unit interval.
	We take the line segment $\mathcal E^\delta_0$ instead and carry out the above Riemann integral construction.
	Evidently, if $v_0=0$, the so redefined $\delta$ is a $\Kinf$ function.
\end{rem}

\begin{rem}
	\label{rem_ifeps0}
	If the uniform overshoot condition in \eqref{eqn_genericlim} read
	\begin{equation*}
			\begin{aligned}
				& \forall \state \in \R^n \spc \lim_{\tau \ra \infty} \Phi(\state, \tau) = 0, \\
				& \exists \eps_0 > 0 \spc \forall \eps \ge \eps_0 \spc \exists \delta \geq 0 \spc \nrm{\state} \le \delta \implies \\
				& \mquad{4} \forall \tau \ge 0 \spc \Phi(\state, \tau) \le \eps,
			\end{aligned}		
	\end{equation*}
	then the $C_0$ would have to be set so as to account for the $\eps_0$ offset.
	In practice, if $\Phi$ were to refer to the distance to some goal under a some given policy, existence of an offset $\eps_0$ would mean that the said policy could not render the goal forward-invariant \ie some overshoot would exist even when starting inside the goal set. 
\end{rem}

For uniform stabilization into the goal, it is crucial to have a uniform $\KL$ bound to estimate overshoot.
By \Cref{prp_lim2kappa}, existence of a uniform stabilizer as per \Cref{dfn_stabilizer} would imply:
\begin{equation}
	\label{eqn_klreachprob}
	\begin{aligned}
		& \forall \state_0 \in \states \\
		& \PP{ \exists \beta \! \in \! \KL \ \goaldist(\State_t) \! \le \! \beta(\goaldist(\state_0)+C_0, t) \mid \Action_t \sim \policy_0(\bullet \mid \state_t)} \\
		& \mquad{10} \geq 1 - \eta,
	\end{aligned}
\end{equation}
where a constant $C_0$ accounts for the $\eps_0$ offset above.

Notice each such $\beta$ depends on the sample variable $\omega$ in general.
This fact would prevent us from uniformly bounding an overshoot.
The reason is that general $\KL$ functions may be rather bizarre despite the common intuition of them being akin to, say, exponentials or something similarly behaved.
In fact, a $\KL$ may entail extremely slow convergence which cannot be described by any polynomial, elementary or even computable functions.
It may formally entail \eg a construction based on the busy beaver function.
We have to avoid such exotic cases.
Hence, we consider special classes of $\KL$ functions which we will refer to as tractable.

First, we need appropriate definitions.
Let us denote a function $\beta$ that depends continuously on some $M \in \N$ parameters besides its main arguments, say, $v, \tau \in \R$ as $\beta\left[\vartheta^{0:M-1}\right](v, \tau)$.
Notice the upper index of the parameter variable $\vartheta$ refers to the respective parameter component, not an exponent. 

\begin{dfn}
	\label{dfn_kappaellfinitely}
	Let $\kappa, \xi$ be of space $\Kinf$ and $M \in \N$.
	A subspace $\KL[\kappa, \xi, M]$ of $\KL$ functions is called \textit{fixed-asymptotic tractable} if it is continuously finitely parametrizable \wrt $\kappa, \xi$ as follows: for any $\beta$ from $\KL[\kappa, \xi, M]$ there exist $M$ real numbers $\vartheta^{0:M-1}$ \sut
	\[
		\forall v, t \spc \beta(v, t) = \kappa\left[\vartheta^{0:M-2}\right](v) \xi\left( e^{-\vartheta^{M-1} t} \right),
	\]
	where $\kappa$ depends continuously on $\vartheta^{0:M-2}$.
\end{dfn}

\begin{exm}
	\label{exm_kappaellfinitely}
	Exponential functions $C v e^{ - \lambda t }$ with parameters $C \ge 0, \lambda \ge 0$ form a fixed-asymptotic tractable $\KL$ subspace.
\end{exm}

The next proposition allows to extract $\omega$-uniform $\KL$ bounds from policies that yield fixed-asymptotic tractable $\KL$ bounds.

\begin{prp}
	\label{prp_fixasymtractpolicy}
	Let $\KL[\kappa, \xi, M]$ be fixed-asymptotic tractable.
	Let $\policy_0$ satisfy

	\begin{equation}
		\label{eqn_fixasymtractstab}
		\begin{aligned}
			\forall \state_0 \in \states \ \mathbb P \big[& \exists \beta \in \KL[\kappa, \xi, M] \spc \goaldist(\State_t) \le \\ & \beta(\goaldist(\state_0)+C_0, t) \mid \Action_t \sim \policy_0(\bullet \mid \state_t) \big] \geq 1 - \eta, \\	
		\end{aligned}
	\end{equation}
	
	where $\eta \in [0,1), \ C_0>0$.Then, for each $\eta'>0$, there exists a $\KL$ function $\beta'$ \sut
	\begin{equation}
		\begin{aligned}	
			\forall \state_0 \in \states \ \mathbb P \big[ & \goaldist(\State_t) \le \beta'(\goaldist(\state_0)+C'_0, t) \\
			&  \mid \Action_t \sim \policy_0(\bullet \mid \state_t) \big] \geq 1 - \eta - \eta', 
		\end{aligned}
	\end{equation}
	where $C'_0>0$.
	Notice that the existence quantifier for $\beta'$ in the latter statement stays before the probability operator, whence uniformity over compacts of the sample space $\Omega$.
\end{prp}

\begin{prf}
	Recall $(\Omega, \Sigma, \mathbb P)$, the probability space underlying \eqref{eqn_mdp} and denote the state trajectory induced by the policy $\policy_0$ emanating from $\state_0$ as $\Traj_{0:\infty}^{\policy_0}(\state_0)$ with single elements thereof denoted $\Traj_t^{\policy_0}(\state_0)$.
	Let us explicitly specify the sample variable in the trajectory $\Traj_{0:\infty}^{\policy_0}(\state_0)$ emanating from $\state_0$ as follows: $\Traj_{0:\infty}^{\policy_0}(\state_0)[\omega], \omega \in \Omega$.
	Define, the subset $\Omega_0$ of elements $\omega \in \Omega$ to satisfy:
	\begin{equation}
		\label{eqn_omegastar}
		\begin{aligned}
			& \Omega_0(\state_0) := \\
			& \left\{ \exists \beta \in \KL[\kappa, \xi, M] \spc \goaldist\left( \Traj_t^{\policy_0}(\state_0)[\omega] \right) \le \beta(\state_0, t) \right\}.
		\end{aligned}
	\end{equation}	
	We argue here assuming a fixed $\state_0$, so we omit the $\state_0$ argument in $\Omega_0$ for brevity.
	In other words, for every $\omega \in \Omega_0$, there exists a $\KL$ function $\beta^\omega$ from $\KL[\kappa, \xi, M]$ \sut 
	\begin{equation}
		\label{eqn_betaomega}
		\forall t \in \Z_{\ge 0} \spc \goaldist\left( \Traj_t^{\policy_0}(\state_0)[\omega] \right) \le \beta^\omega(\goaldist(\state_0), t).
	\end{equation}
	By the condition of the lemma, $\PP{\Omega_0} \ge 1 - \eta$.

	Now, let $\vartheta_{k=0:\infty}$ be a dense sequence in $\R^{M}$.
	
	Then, let $\{\beta_k\}_{k=0:\infty}$ be the sequence defined by $\beta_k(v, t) := \kappa\left[\vartheta^{0:M-2}_k\right](v) \xi\left( e^{-\vartheta^{M-1}_k t}, \right)$ for all $v, t$.
	First, observe that $\{\beta_k\}_{k=0:\infty}$ is dense in $\image(\beta[ \bullet ])$ since the dependence of $\kappa\left[\vartheta^{0:M-2}_k\right](v) \xi\left( e^{-\vartheta^{M-1}_k t} \right)$ on the parameters $\vartheta^{0:M-1}$ is continuous and continuous functions map dense sets into dense sets.
	
	Define a sequence of functions $\{\hat \beta_k\}_{k=0:\infty}$ by $\hat \beta_k(v, t) := \left( \kappa\left[\vartheta^{0:M-2}_k\right](v) + 1 \right) \xi\left( e^{-\vartheta^{M-1}_k t} \right)$ for all $v, t$.
	These are $\KL$ functions with an offset along the $v$-dimension.

	The reason for such an offset will be revealed in \eqref{eqn_supbound}.
	
	
	Let $\hat \Omega_0^{(k)}, j \in \Z_{\ge 0}$ be defined via
	\[
	\forall \omega \in \hat \Omega_0^{(k)} \spc \forall \state_0 \in \states, t \ge 0 \spc \goaldist\left( \Traj_t^{\policy_0}(\state_0)[\omega] \right) \le \hat \beta_k(\goaldist(\state_0), t).
	\]
	Then, 
	\[
		\hat \Omega_0^{(\infty)} := \bigcap_{k=0:\infty} \hat \Omega_0^{(k)}
	\]
	is an event, since it is a countable intersection.
	Now, to any $\beta^\omega$ there corresponds some $\hat \beta_k$ \sut $\hat \beta_k \ge \beta^\omega$.
	To see this, consider
	\[
		\beta^\omega(v, t) = \kappa\left[\vartheta^{0:M-2}\right](v) \xi\left( e^{-\vartheta^{M-1} t} \right)
	\] 
	as a function of $v, t$ with fixed $\vartheta^{0:M-1}$.
	By the density of $\{\beta_k\}_{k=0:\infty}$ in $\image(\beta[ \bullet ])$, there exists some $\beta_l, l \in \Z_{\ge 0}$ \sut 

	\begin{equation}
		\label{eqn_supbound}
		\begin{aligned}
				\sup_{v, t} \bigg\| & \kappa\left[\vartheta^{0:M-2}\right](v) \xi\left( e^{-\vartheta^{M-1} t} \right) - \\
				& \kappa\left[\vartheta^{0:M-2}_l\right](v) \xi\left( e^{-\vartheta^{M-1}_l t} \right) \bigg\| \le \frac{\xi(1)}{2}.
		\end{aligned}
	\end{equation}

	Now, take $\vartheta_k$ \sut $\vartheta^{0:M-2}_k := \vartheta^{0:M-2}_l$ and $\vartheta^{M-1}_k \le \vartheta^{M-1}$.
	From \eqref{eqn_supbound}, it is now clear how an offset by one in $\hat \beta$s along with the approximation bound $\frac{\xi(1)} {2}$ ensure $\hat \beta_k \ge \beta^\omega$.
	Namely, \eqref{eqn_supbound} implies
	\[
		\sup_{v} \nrm{ \kappa\left[\vartheta^{0:M-2}\right](v) - \kappa\left[\vartheta^{0:M-2}_l\right](v) } \le \frac{1}{2}.
	\]
	This in turn implies that
	\[
		\forall v \spc \kappa\left[\vartheta^{0:M-2}\right](v) \le \kappa\left[\vartheta^{0:M-2}_l\right](v) + 1.
	\]
	This together with the fact that
	\[
		\forall t \spc \xi\left( e^{-\vartheta^{M-1} t} \right) \le \xi\left( e^{-\vartheta_k^{M-1} t} \right)
	\]
	ensures $\hat \beta_k \ge \beta^\omega$.
	
	
	

	Then, it holds that $\Omega_0 \subseteq \hat \Omega_0^{(\infty)}$.
	Define a new sequence $\{\bar \beta_k\}_{k=0:\infty}$ as follows:
	\begin{equation}
		\label{eqn_barbetasqn}
		\begin{aligned}
			\forall v, t \spc \bar \beta_k(v, t) := \sum_{j \le k} \hat \beta_j(v, t),
		\end{aligned}
	\end{equation}	
	which is still a $\KL$ function with an offset along the $v$-dimension since it is a sum of finitely many $\KL$ functions with offsets along the $v$-dimension.
%
%
	Then it holds that $\bar \beta_{k+1} \ge \bar \beta_k$ meaning $\forall v \in \R_{\ge 0}, \forall t \ge 0 \spc \bar \beta_{k+1}(v, t) \ge \bar \beta_{k}(v, t)$.
	Then, evidently
	\begin{equation}
		\label{eqn_omegastarin}
		\Omega_0 \subseteq \bar \Omega_0^{(\infty)} := \bigcap_{k=0:\infty} \bar \Omega_0^{(k)},
	\end{equation}
	where $\bar \Omega_0^{(k)}$ is defined similarly to $\hat \Omega_0^{(k)}$ with $\bar \beta_k$ in place of $\hat \beta_k$.
	
	Too see this, fix an arbitrary $k \in \Z_{\ge 0}$ and observe that for any $\omega$
	\[
		\forall t \in \Z_{\ge 0}, \state_0 \in \states \spc \goaldist\left( \Traj_t^{\policy_0}(\state_0)[\omega] \right) \le \bar \beta_k(\goaldist(\state_0))
	\]
	implies
	\[
		\forall t \in \Z_{\ge 0}, \state_0 \in \states \spc \goaldist\left( \Traj_t^{\policy_0}(\state_0)[\omega] \right) \le \hat \beta_k(\goaldist(\state_0)),
	\]
	hence an event $\bar \Omega_0^{(k)}$ implies an event $\hat \Omega_0^{(k)}$.
	
	By the construction \eqref{eqn_barbetasqn}, 
	\[
		\forall k \in \Z_{\ge 0} \spc \bar \Omega_0^{(k)}  \subseteq \bar \Omega_{k+1}.
	\]
	Furthermore, by \eqref{eqn_omegastarin}
	\begin{equation}
		\label{eqn_limbaromegaprob}
	   \PP{\bar \Omega_\infty} \ge \PP{\Omega_0}. 
	\end{equation}
	
	Notice $\bar \beta_k$ was constructed uniform on the respective $\bar \Omega_0^{(k)}$ in terms of $\omega$.
	Intuitively, $\bar \Omega_0^{(k)}$ are the sets of outcomes for which there is a uniform convergence modulus.
	
	For any $k \in \Z_{\ge 0}$ we have $\PP{\bar \Omega_0^{(k)}} \ge 1 - \eta + \eta_k$ where $\eta_k$ is a real number.

	The condition $\PP{\bar \Omega_0^{(k)}} \le \PP{\bar \Omega_{k+1}}$ implies that the sequence $\{\eta_k\}_{k=0:\infty}$ is a monotone increasing sequence which is bounded above.
	Hence, $\lim_{k \ra \infty} \eta_k$ exists.
	It also holds that $\lim_{k \ra \infty} \eta_k \le 0$ by \eqref{eqn_limbaromegaprob}.

	We thus argue that with probability not less than $1 - \eta + \eta_k$, the policy $\policy_0$ is a uniform stabilizer with an $\omega$-uniform $\KL$ bound.
	
	Fix any $k \in \Z_{\ge 0}$ and take the corresponding $\bar \beta_k$.
	It has the form
	\[
		\bar \beta_k(v, t) = \sum_{j \le k} \left( \kappa\left[\vartheta^{0:M-2}_j\right](v) + 1 \right) \xi\left( e^{-\vartheta^{M-1}_j t} \right).
	\]
	Let $\tilde \beta_k$ be defined as
	\[
		\tilde \beta_k(v, t) = \sum_{j \le k} \left( \kappa\left[\vartheta^{0:M-2}_j\right](v+C'_0) \right) \xi\left( e^{-\vartheta^{M-1}_j t} \right),
	\]
	where $C'_0$ is chosen \sut $\tilde \beta_k \equiv \bar \beta_k$.
	
	Finally, $\tilde \beta$ is a $\KL$ function sought to satisfy the claim of the proposition.
	
\end{prf}

\begin{rem}
	\label{rem_zerooffset}
	Let us assume $C_0=0$.
	The construction $\hat \beta_k(v, t) := \left( \kappa\left[\vartheta^{0:M-2}_k\right](v) + 1 \right) \xi\left( e^{-\vartheta^{M-1}_k t} \right)$ for all $v, t$, necessarily leads to a non-zero offset $C_0'$ in the final $\KL$ function $\tilde \beta$, namely, it satisfies
	\[
		\sum_{j \le k} \kappa\left[\vartheta^{0:M-2}_j\right](C'_0) \xi\left( 1 \right) = k \xi(1).
	\]
	We cannot in general eliminate such an offset -- the density argument does not suffice here since for any fixed sequence of $\KL$ functions, there can always be one that grows faster than any selected one from the sequence in a vicinity of zero.
	However, the offset can be made arbitrarily small leading to an arbitrary small relaxation of the goal set to which stabilization is of interest.
	
	This can be done \eg by setting
	\[
		\hat \beta_k(v, t) := \left( \kappa\left[\vartheta^{0:M-2}_k\right](v) + \lambda_\kappa^{k+1} \right) \xi\left( e^{-\vartheta^{M-1}_k t} \right),
	\]
	where $0<\lambda_\kappa<1$ can be chosen arbitrarily small.
	Then, the offset will not exceed
	\[
		\frac{\lambda_\kappa}{1-\lambda_\kappa} \xi(1),
	\]
	leading to the respective arbitrarily small constant $C'_0$.
\end{rem}

Now, we address subspaces of $\KL$ functions called \textit{weakly tractable}.
First, let introduce the following relation on the space of two scalar argument functions:

\begin{dfn}
	\label{dfn_orderonkappas}
	A function of two arguments $\beta_2$ is said to precede another function of two arguments $\beta_1$, with the respective relation denoted as $\beta_1 \succ \beta_2$, if $\forall v, t \spc \beta_1(v,t) > \beta_2(v, t)$.
\end{dfn}

Next, we consider a special notion of continuity on the space of parametric maps into $\KL$ functions made handy afterwards.
Let $\mathcal P = \R^M, M \in \N$.

\begin{dfn}
	\label{dfn_klcontparammaps}
	A parametric map $\beta[\bullet]: \mathcal P \ra \KL$ is called $\KL$-continuous if for any $\KL$ function $\xi$ there exists $\delta>0$ \sut $\forall \vartheta_1, \vartheta_2 \in \mathcal P$ with $\nrm{\vartheta_1 - \vartheta_2} \le \delta$ it holds that $\abs{ \beta[\vartheta'] - \beta[\vartheta] } \prec \xi$.
\end{dfn}

Finally, we introduce proper $\KL$-continuous maps which avoid pathological cases where $\beta[\bullet]$ could get ``saturated'' with the parameter, in other words, have a finite limit superior in the $\KL$-continuity sense.

\begin{dfn}
	\label{dfn_klproperparammaps}
	A parametric map $\beta[\bullet]: \mathcal P \ra \KL$ is called proper if it is continuous and for any $\vartheta \in \mathcal P$ there exists a $\vartheta' \in \mathcal P$ \sut $\beta[\vartheta'] \succ \beta[\vartheta]$.
\end{dfn}

Notice that the condition $\beta[\vartheta'] \succ \beta[\vartheta]$ implies existence of a $\KL$ function $\xi$ \sut $\beta[\vartheta'] \succeq \beta[\vartheta]+\xi$.
Take \eg $\xi := \frac{\beta[\vartheta'] + \beta[\vartheta]}{2}$.

\begin{dfn}
	\label{dfn_weaktractkl}
	A subspace $\mathcal F$ of $\KL$ is called weakly tractable if there exists a proper parametric map $\beta[\bullet]: \mathcal P \ra \KL$ \sut $\image(\beta[\bullet]) = \mathcal F$.
\end{dfn}

\begin{exm}
	A fixed-asymptotic tractable $\KL$ subspace consisting of functions of the form
	\[
		\forall v, t \spc \beta(v, t) = \kappa\left[\vartheta^{0:M-2}\right](v) \xi\left( e^{-\vartheta^{M-1} t} \right),
	\]	
	where $\xi$ is unbounded \wrt $\vartheta^{M-1}$, and for any $\vartheta_1$ there exists $\vartheta_2$ \sut $\forall v \spc \kappa\left[\vartheta_1^{0:M-2}\right](v) < \kappa\left[\vartheta_2^{0:M-2}\right](v)$, is weakly tractable.
	However, a general fixed-asymptotic tractable $\KL$ subspace is not necessarily weakly tractable.
\end{exm}

Now, the following may be claimed.

\begin{prp}
	\label{prp_weaktractpolicy}
	Let $\mathcal F$ be weakly tractable.
	Let $\policy_0$ satisfy
	\begin{equation}
		\label{eqn_weaktractstab}
		\begin{aligned}
			& \forall \state_0 \in \states \\
			& \PP{ \exists \beta \in \mathcal F \spc \goaldist(\State_t) \le \beta(\goaldist(\state_0)+C_0, t) \mid \Action_t \sim \policy_0(\bullet \mid \state_t)} \\
			& \pushright{\geq 1 - \eta, \eta \in (0,1), C_0>0.}	
		\end{aligned}
	\end{equation}
	
	Then, for each $\eta'>0$, there exists a $\KL$ function $\beta'$ \sut
	\begin{equation}
		\begin{aligned}
			& \forall \state_0 \in \states \\
			& \PP{\goaldist(\State_t) \le \beta'(\goaldist(\state_0)+C_0, t) \mid \Action_t \sim \policy_0(\bullet \mid \state_t)} \\
			& \pushright{\geq 1 - \eta - \eta'.} 		
		\end{aligned}
	\end{equation}
	Notice that the existence quantifier for $\beta'$ in the latter statement stays before the probability operator, whence uniformity over compacts of the sample space $\Omega$.
\end{prp}

\begin{prf}
	The proof follows the same ideas as the one for \Cref{prp_fixasymtractpolicy}, yet many constructions become simpler.
	We again choose a dense sequence $\{\vartheta_k\}_k$ in $\mathcal P$.
	Then, for a given $\beta^\omega$ from $\mathcal F$, we find a $\hat \beta \in \mathcal F$ \sut $\hat \beta \succeq \beta^\omega + \xi$ for some $\KL$ function $\xi$ using the fact that $\beta[\bullet]$ is proper. 
	Using $\KL$-continuity of $\beta[\bullet]$, take a $\xi/2$-close neighbor to $\hat \beta$ from $\{\beta[\vartheta_k]\}_k$.
	The rest of the proof is essentially the same as for \Cref{prp_fixasymtractpolicy}.	
\end{prf}

\begin{rem}
	\label{rem_zerooffsetweaktract}
	Notice that the final $\KL$ function possesses the same offset as in the condition \eqref{eqn_weaktractstab}, in contrast to the fixed-asymptotic case.
\end{rem}

Let us refer to policies which yield fixed-asymptotic or weakly tractable $\KL$ bounds as per \eqref{eqn_fixasymtractstab}, \eqref{eqn_weaktractstab} fixed-asymptotic and, respectively, weakly tractable stabilizers.

\begin{exm}
	\label{exm_expstab}
	For a controlled Markov chain 
	\begin{equation}
		\label{eqn_linmdp}
		\begin{aligned}
			& \State_{t + 1} = F \State_t + G \Action_t + \Noise_t, \\
			& \qquad F, G \text{ being matrices of proper dimension},
		\end{aligned}
	\end{equation}
	where $\Noise_t$ is uniformly distributed on $[-\bar \noise, \bar \noise]^n, n \in \N, \bar \noise < \diam_\G$, a condition
	\begin{equation}
		\label{eqn_expstab}
		\PP{\lim_{t \ra \infty} \goaldist(\State_t) = 0 \mid \Action_t \sim \policy_0(\bullet \mid \state_t)} \ge 1-\eta, \eta \in (0,1)
	\end{equation}	
	implies $\policy_0$ is a stabilizer with the respective exponential fixed-asymptotic tractable $\KL$ subspace as per \Cref{exm_kappaellfinitely}.
\end{exm}

Note that known converse results cannot, in general, be directly applied to the above setting, due to arbitrariness of $\mathbb{G}$ and also the possibility of failed stabilization.

We are in the position to prove the following theorem, seen as a stochastic generalization of \cite{Jiang2002converseLyapun,Sontag1998Commentsintegr}.

\begin{thm}[Tractable stochastic converse LF]
	\label{thm_problfconverse}
	Let $\policy_0 \in \policies_0$ be an $\eta$-improbable fixed-asymptotic or weakly tractable $\G$-stabilizer for the controlled Markov chain \eqref{eqn_mdp}.
	Then, for any $\eta'>0$ there exist a locally Lipschitz-continuous function $L:\R^n \ra \R_{\ge 0}$, a continuous function $\nu:\R^n \ra \R_{\ge 0}$ that is strictly positive outside of $\G$, and $\kappa_\low, \kappa_\up \in \Kinf$ \sut 
	\begin{enumerate}
		\item[i)] $L$ has a decay property satisfying
			\begin{align*}
					& \forall \state_t \notin \G', \Action_t \sim \policy(\bullet \mid \State_t) \implies \\
					& \quad \PP{L(\State_{t + 1}) - L(\state_{t}) \le -\nu(\goaldist['](\state))} \ge 1 -\eta-\eta',
			\end{align*}
			where $\G \subseteq \G'$,
		\item[ii)] $\forall \state \in \states \ \ \kappa_\low(\goaldist['](\state)) \leq L(\state) \le \kappa_\up(\goaldist['](\state))$.
	\end{enumerate}
\end{thm}

\begin{prf}
	Let $(\Omega, \Sigma, \PP{\bullet})$ be the probability space underlying \eqref{eqn_mdp} and denote the state trajectory induced by the policy $\policy$ emanating from $\state_0$ as $\Traj_{0:\infty}^{\policy}(\state_0)$ with single elements thereof denoted $\Traj_t^{\policy}(\state_0)$.
	
	By \Cref{prp_fixasymtractpolicy} or \Cref{prp_weaktractpolicy} (depending on fixed-asymptotic or weak tractability of $\policy_0$), it holds that there exists a $\KL$ function $\beta$ satisfying:
	\begin{equation}
		\label{eqn_klboundinthm}
		\begin{aligned}
			& \forall \state_0 \in \states \\
			& \PP{\goaldist(\State_t) \le \beta(\goaldist(\state_0)+C'_0, t) \mid \Action_t \sim \policy_0(\bullet \mid \state_t)} \\
			& \pushright{\ge 1 - \eta - \eta', C'_0>0,}
		\end{aligned}
	\end{equation}		
	where $C'_0 = C_0$ if $\policy_0$ is weakly tractable.
	Let $C'_0, \G'$ be \sut $\beta(\goaldist(\state_0)+C'_0, t) = \beta(\goaldist['](\state_0), t)$ for all $\state_0, t$.
	In particular, $\G'=\G$ if $C_0=0$ and $\policy_0$ is weakly tractable.

	by Proposition 7 in \cite{Sontag1998Commentsintegr}, for $\beta$ there exist $\kappa_1, \kappa_2 \in \Kinf$ \sut
	\begin{equation}
		\label{eqn_klcompose}
		\beta(v, t) = \kappa_1( \kappa_2 (v)e^{-t}) \spc \forall v > 0,  t > 0.
	\end{equation}
	
	Let $\cost' := \kappa_1\inv$. Using \eqref{eqn_klboundinthm} and \eqref{eqn_klcompose}, we obtain
		\begin{align*}
			& \PP{\goaldist['](\Traj_t^{\policy_0}(\state_0)) \le \kappa_1 \left( \kappa_2 \left(\goaldist['](\state_0) \right) e^{-t+1} \right) } \ge 1 - \eta - \eta',
		\end{align*}
	which in turn implies
	\begin{equation}
		\label{eqn_expbound}
		\PP{ \cost'(\goaldist['](\Traj_t^{\policy}(\state_0))) \le \left( \kappa_2 \left(\goaldist['](\state_0) \right) \right) e^{-t+1} } \ge 1 - \eta - \eta'.
	\end{equation}
	Define $L' : \states \times \states^\infty \ra \R_{\ge 0}$ by
	\begin{equation}
		\label{eqn_lfseries}
		L' [\traj_{0:\infty}(\state_0)] := \sum_{t=0}^{\infty}\cost'(\goaldist['](\traj_t(\state_0))),
	\end{equation}
	with an emphasis on the initial state $\state_0$.
	Denote $\kappa_\low(\goaldist['](\state_0)) := \cost'(\goaldist['](\state_0))$ and $\kappa_\up(\goaldist['](\state_0)) := \frac{e^2}{e - 1}\kappa_2(\goaldist['](\state_0))$.
	Then, from \eqref{eqn_expbound} it follows that, with probability no less than $1 - \eta - \eta'$:
	
	\begin{equation*}
	\begin{aligned}
		& \kappa_\low(\goaldist['](\state_0)) = \cost'(\goaldist['](\state_0)) \le L'[\Traj_{0:\infty}^{\policy}(\state_0)] \le \\
		& \sum_{t=0}^{\infty}\!\kappa_2(\goaldist['](\state_0))e^{-t+1} \!\le\! \frac{e^2}{e - 1}\kappa_2(\goaldist['](\state_0)) \!=\! \kappa_\up(\goaldist['](\state_0)).
	\end{aligned}
	\end{equation*}	
	
	The above also shows that the series in \eqref{eqn_lfseries} converges uniformly on compact sets with respect to $\state_{0}$, because the tail sums can be bounded in the following way ($T$ is arbitrary) with probability no less than $1 - \eta - \eta'$:
	\begin{equation}
		\label{eqn_lfseriesconv}
		\begin{aligned}
			& 0 \le \sum_{t=T}^{\infty} \cost'(\goaldist['](\Traj_t^{\policy})) \le \sum_{t=T}^{\infty}\sup\limits_{\state \in \states_0} \kappa_2(\goaldist['](\state))e^{-t+1} \le \\
			& \mquad{1} \frac{e}{1 - e^{-1}}\sup\limits_{\state \in \states_0}\kappa_2(\goaldist['](\state)) - \frac{e - e^{1-T}}{1 - e^{-1}}\sup\limits_{\state \in \states_0}\kappa_2(\goaldist['](\state)) \\
			& \mquad{9} \xrightarrow{T \ra \infty} 0,
		\end{aligned}
	\end{equation}
	where $\states_0$ is an arbitrary compact set containing $s_0$.	
	Thus, the series in \eqref{eqn_lfseries} converges uniformly with respect to $\state_t, t \in [0,\infty)$, while the elements of the sequence are Lipschitz continuous with respect to $\state_0$.
	Thus, by the uniform convergence theorem the series converges to a function that is uniformly continuous with respect to $\state_0$.
	Let us briefly and explicitly specify the sample variable in a trajectory $\Traj_{0:\infty}^{\policy}(\state)$ emanating from $\state$ as follows: $\Traj_{0:\infty}^{\policy}(\state)[\omega], \omega \in \Omega$.
	Let $\Omega_0'(\state)$ be the event corresponding to \eqref{eqn_klboundinthm} with the starting state $\state$.
	Define $L(\state) := \esssup_{\omega \in \Omega_0'(\state)} L' [\Traj_{0:\infty}^{\policy}(\state)[\omega]]$ and
	\begin{equation}
		\begin{aligned}
			\kappa_\low:= \cost', \quad \kappa_\up := \frac{e^2}{e - 1}\kappa^{\#}_2, \quad \nu := \cost'.
		\end{aligned}
	\end{equation}		
	These functions constitute what was sought.	
\end{prf}


Let ``$\bullet \Subset \bullet$'' denote a subset inclusion with a non-zero Hausdorf distance between the included and container sets.

\begin{prp}[Mean LF]
	\label{lem_meanlfconverse}
	Let $\policy_0$ be a uniform mean stabilizer in the following sense:

	\begin{equation}
		\label{eqn_meanstabilizer}
		\exists \beta \in \KL \spc \E[\policy_0]{\goaldist(\State_t)} \le \beta(\nrm{\state_0}, t) \ \forall \state_0 \in \states, t \ge 0.
	\end{equation}

	Then, for any $\G' \Subset \G$ there exist a locally Lipschitz-continuous function $L:\R^n \ra \R_{\ge 0}$, a continuous function $\nu:\R^n \ra \R_{\ge 0}$ that is strictly positive outside of $\G'$, and $\kappa_\low, \kappa_\up \in \Kinf$ \sut 
	\begin{enumerate}
		\item[i)] $L$ has a decay rate satisfying
			\begin{align*}
					& \forall \state_t \notin \G', \Action_t \sim \policy(\bullet \mid \state_t) \implies \\
					& \quad \E{L(\State_{t + 1})} - L(\state_{t}) \le -\nu(\goaldist['](\state)),
			\end{align*}
		\item[ii)] $\forall \state \in \states \ \ \kappa_\low(\goaldist['](\state)) \leq L(\state) \leq \kappa_\up(\goaldist['](\state))$.
	\end{enumerate}
\end{prp}

\begin{prf}
	Denote, similar to Theorem \ref{thm_problfconverse}, the state trajectory induced by the policy $\policy$ emanating from $\state_0$ as $\Traj_{0:\infty}^{\policy}(\state_0)$ with single elements thereof denoted $\Traj_t^{\policy}(\state_0)$.	
	
	By Proposition 7 in \cite{Sontag1998Commentsintegr}, for the $\KL$ function $\beta_E$ there exist $\kappa_1, \kappa_2 \in \Kinf$ \sut
	\begin{equation}
		\label{eqn_meanklcompose}
		\beta_E(v, t) = \kappa_1( \kappa_2 (v)e^{-t}) \spc \forall v > 0,  t > 0.
	\end{equation}
	These $\kappa_1, \kappa_2$ are standalone here and not to be confused with the ones in Theorem \ref{thm_problfconverse}.
	Now, since $\kappa_1$ is $\Kinf$, for an arbitrary $\bar v > 0$ there exists another $\Kinf$-function $\kappa_1'$, which is convex, and it holds that $\forall v \in [\bar v, \infty) \spc \kappa_1(v) \le \kappa_1'(v)$.
	Let $\eps_\G$ be the Hausdorf distance between $\G$ and $\G'$, and choose $\bar v > 0 $ \sut the corresponding $\kappa_1'$ yields
	\[
		\forall \state \in \states, t>0 \spc \beta_E(\goaldist['](\state), t) \le \kappa_1'( \kappa_2 (\goaldist['](\state))e^{-t} ).
	\]
	
	Let $\cost' := \left(\kappa_1'\right)\inv$.
	Notice $\cost'$ is concave because $\kappa_1'$ is convex and strictly increasing.
	
	Under \eqref{eqn_meanstabilizer} and \eqref{eqn_meanklcompose}, it holds that
	\begin{align*}
		& \E[\policy_0]{\goaldist['](\Traj_t^{\policy}(\state_0))} \le \beta_E(\goaldist['](\state_0), t-1) \le \\
		& \mquad{5} \kappa'_1(\kappa_2(\goaldist['](\state_0))e^{-t+1}).
	\end{align*}
	This implies
	\begin{equation}
		\label{eqn_meanexpbound}
		\cost'\left( \E[\policy_0]{ \goaldist['](\Traj_t^{\policy}(\state_0)) } \right) \le \kappa_2(\goaldist['](\state_0))e^{-t+1}.
	\end{equation}
	By the Jensen's inequality, we have (noticing that $\cost'$ is concave):
	\begin{equation}
		\label{eqn_Jensen}
		\E[\policy_0]{ \cost'\left( \goaldist['](\Traj_t^{\policy}(\state_0)) \right) } \le \cost'\left( \E[\policy_0]{ \goaldist['](\Traj_t^{\policy}(\state_0)) } \right).
	\end{equation}	
	
	Define $L: \states \ra \R_{\ge 0}$ by $L(\state) := \sum_{t=0}^{\infty}\E[\policy_0]{ \cost'\left( \goaldist['](\Traj_t(\state)) \right)}$.	
	It holds that, defining $\kappa_\low := \cost'$,
	\begin{equation}
		\label{eqn_meankappalow}
		L(\state) \ge \E[\policy_0]{ \cost'\left( \goaldist[']( \state \right) } = \cost'\left( \goaldist[']( \state ) \right) = \kappa_\low \left( \goaldist[']( \state ) \right).
	\end{equation}
	On the other hand, from \eqref{eqn_meanexpbound}, we have
	\begin{equation}
		\label{eqn_meankappaup}
		L(\state) \le \sum_{t=0}^{\infty} \kappa_2(\goaldist(\state))e^{-t+1} \le \frac{e^2}{e - 1}\kappa_2(\goaldist(\state)).
	\end{equation}
	The uniform convergence follows analogously to \eqref{eqn_lfseriesconv} taking $\sum_{t=T}^{\infty} \E[\policy_0]{ \cost'\left( \goaldist['](\Traj_t(\state)) \right) }$ as the second term in the respective chain of inequalities.
	Setting $\kappa_\up := \frac{e^2}{e - 1}\kappa^{\#}_2, \nu := \cost'$ concludes the proof noticing that

	\begin{multline}
			\E[\policy_0]{L(\State_1)} - L(\state_0) = \\
			 \sum_{t=1}^{\infty}\E[\policy_0]{ \cost'\left( \goaldist['](\Traj_t(\state_0)) \right)} - \sum_{t=0}^{\infty}\E[\policy_0]{ \cost'\left( \goaldist['](\Traj_t(\state_0)) \right)} = \\
		  \E[\policy_0]{ \cost'\left( \goaldist['](\Traj_0(\state_0) \right) } = -\cost'(\goaldist['](\state_0)). 
		\end{multline}

\end{prf}


\subsection{Corresponding direct Lyapunov results}

In this section, we briefly discuss the opposite direction, namely, given a stochastic Lyapunov function, derive a stabilizer (not knowing the original stabilizer which the Lyapunov function was constructed from).
Let us first assume a control Lyapunov function satisfying a decay similar to that of Theorem \ref{thm_problfconverse}, namely, omitting $\eta'$ and assuming $\G=\G'$ for brevity,
\begin{equation}
	\label{eqn_espclf}
	\begin{aligned}
		& \forall \state \in \states \spc \exists \action \\
		& \PP{L(\State_+) - L(\state) \!\le\! - \nu(\goaldist(\state)) \mid \State_+ \sim \transit(\bullet \mid \state, \action)} \!\ge\! 1-\eta.
	\end{aligned}
\end{equation}

Let us assume $\G$ to be a compact vicinity of the origin.
Let us take 
\begin{multline}
	\label{eqn_steepestdescent}
	\action_t \la \max_\action \mathbb{P}\Big[(\State_{t+1}) - L(\state_t) \le \\ - \nu(\goaldist(\state_t)) \mid \State_{t+1} \sim \transit(\bullet \mid \state_t, \action)\Big],
\end{multline}
which is akin to a steepest descent policy (again, our goal here is not to take the fanciest algorithm).
Following this policy, which we denote $\policy'$, we have:
	\begin{equation}
	\forall t \in \Z_{\ge 0} \ \PP{L( \goaldist(\Traj_t^{\policy'}(\state_0)) \! \le \! \kappa_\up( \goaldist(\state_0) ) } \! \ge \! 1 - \eta.
	\end{equation}
Hence, with probability not less than $1-\eta$ the state is uniformly stable \ie
\begin{equation}
	\label{eqn_epsstable}
	\begin{aligned}
	& \forall t \in \Z_{\ge 0} \\
	& \mquad{1} \PP{ \nrm{\Traj_t^{\policy}(\state_0)} \le \kappa_\low\inv \left(\kappa_\up( \nrm{\state_0} )\right) + \diam_\G } \ge 1 - \eta.
	\end{aligned}
\end{equation}
Monotonicity of $\kappa_\low\inv$ was employed here.
Convergence to the goal set with probability not less than $1-\eta$ follows considering an arbitrary superset $\G' \Subset \G$.
Then, if a state $\State_t$ is not in $\G'$, it holds that $\nu(\goaldist(\State_t)) \ge \inf_{\state \notin \G'} \nu(\goaldist(\state))$.
Then, $\G'$ is reached with probability not less than $1-\eta$ in no greater than $T_L := \frac{\kappa_\up( \goaldist(\state_0) ) - \kappa_\low\left(\inf\limits_{\state \notin \G'} \goaldist(\state) \right) }{\inf\limits_{\state \notin \G'} \nu(\goaldist(\state))}$ steps.
This holds since with probability not less than $1-\eta$, there is a decay of $L$ as soon as $\State_t$ is not in $\G'$
But since $\G'$ was an arbitrary superset, we conclude asymptotic convergence to $\G$ with an $\omega$-uniform rate.
The above argumentation is quite similar to the deterministic case, except for the probability requirement relaxation (from strict to non-strict in the construction of the Lyapunov function).
The crucial step is \eqref{eqn_epsstable}, where we relied on monotonicity of $\kappa_\low\inv$ since it is an inverse of a $\K$ function.
The situation becomes much more difficult when we consider stabilization in mean using a Lyapunov function as constructed in Lemma \ref{lem_meanlfconverse}. 
Notice that to claim
\begin{equation}
	\label{eqn_meanstable}
	\begin{aligned}
	& \forall t \in \Z_{\ge 0} \spc \E[\policy]{ \nrm{\Traj_t^{\policy}(\state_0)} }\! \le\! \kappa_\low\inv \left(\kappa_\up( \goaldist['](\state_0) )\right) \! + \! \diam_{\G'},
	\end{aligned}
\end{equation}
we have to employ Jensen's inequality.
But it would in turn put the requirement for $\kappa_\low\inv$ to be concave.
That would in turn mean $\kappa_\low$ to be convex.
But $\kappa_\low$ in Lemma \ref{lem_meanlfconverse} was constructed to be concave.
Thus, every cycle of constructing a mean Lyapunov function from a mean stabilizer and vice versa ``straightens'' the lower kappa, since only a linear function is both convex and concave.
Such a behavior of mean is not surprising and ``convexification'' is not unusual -- see \eg constructions in \cite{Subbaraman2013converseLyapun}.

Finally, let us remark on the tractability of stabilizers.
If the tractability property were omitted, we believe that no Lyapunov function construction could be possible.
Consider \eg the property that
\begin{equation}
	\label{eqn_nontractablestabilizer}
	\begin{aligned}
		\PP{\exists \beta \in \KL \spc \goaldist(\State_t) \le \beta(\state_0, t) \mid \Action_t \sim \policy_0(\bullet \mid \state_t)} > 1 - \eta.	
	\end{aligned}
\end{equation}
If $\beta$ is an arbitrary $\KL$ function, it is impossible to derive $\omega$-uniform convergence moduli.
Essentially, for any unbounded sequence $\beta_{k=0:\infty}$ of $\KL$ functions which is dense, by the diagonal argument there exists a $\KL$ function that intersects any of $\beta_{k=0:\infty}$.
One cannot in turn rely on sequences indexed by \eg a real since the uncountable intersection of the respective events needs not to be an event.
Hence, although a full formal negative result is beyond the scope of this work, we strongly suggest that tractability is needed to derive Lyapunov functions from stabilizers.

%

\section{Conclusion}

We presented several elementary constructions of stochastic Lyapunov functions following the classical recipe of \cite{Jiang2002converseLyapun}.
It was shown that a stochastic Lyapunov function could be derived for stabilizers with tractable stabilization certificates.
It was argued, although not formally proven, that no effective Lyapunov function construction could be carried out for stabilizers with arbitrary stabilization certificates that depend on the outcome.
It was also discussed that the constructed mean Lyapunov function possessed rather restrictive properties in terms of the ability to derive a stabilizer (the reverse direction from a Lyapunov function to a stabilizer).
%
%


\bibliography{
bib/AIDA__May2024,
bib/Osinenko__May2024
}

\end{document}